\documentclass[12pt]{amsart}
\usepackage{graphicx}
\usepackage{amscd}
\usepackage[top=1in,bottom=1in,left=1.15in,right=1.15in]{geometry}
\usepackage{amsmath,amssymb,amsfonts}
\usepackage{amsthm}
\usepackage{color}
\usepackage{hyperref}
\usepackage{tikz-cd}
\usetikzlibrary{plotmarks}
\usepackage{epstopdf}
\usetikzlibrary{matrix}
\usepackage{verbatim}
\newtheorem{theo}{Theorem}[section]
\newtheorem{prop}[theo]{Proposition}
\newtheorem{lemma}[theo]{Lemma}
\newtheorem{defn}[theo]{Definition}
\newtheorem{rem}[theo]{Remark}

\newtheorem{ex}[theo]{Example}

\newcommand{\cal}{\mathcal}

\allowdisplaybreaks

\def\diaCrossP{\unitlength.08em
  \begin{minipage}{15\unitlength}
    \begin{picture}(15,15)
      \put(0,0){\vector(1,1){15}}
      \qbezier(15,0)(15,0)(10,5)
      \qbezier(5,10)(0,15)(0,15)
      \put(0,15){\vector(-1,1){0}}
    \end{picture}
  \end{minipage}
}

\def\diaCrossN{\unitlength.08em
  \begin{minipage}{15\unitlength}
    \begin{picture}(15,15)
      \put(15,0){\vector(-1,1){15}}
      \qbezier(0,0)(0,0)(5,5)
      \qbezier(10,10)(15,15)(15,15)
      \put(15,15){\vector(1,1){0}}
    \end{picture}
  \end{minipage}
}

\usetikzlibrary{decorations.markings}

\tikzset{->-/.style={decoration={
  markings,
  mark=at position .5 with {\arrow{>}}},postaction={decorate}}}

\tikzset{-<-/.style={decoration={
  markings,
  mark=at position .5 with {\arrow{<}}},postaction={decorate}}}
  
\def\lgraph{\unitlength.1em
  \begin{minipage}{25\unitlength}
    \begin{picture}(25,19)
      \put(5,2){\line(1,0){11}}
      \put(5,2){\line(0,1){11}}
      \put(5,2){\line(5,2){11}}
      \put(18,2){\textbf {\tiny 1}}
      \put(17,7){\textbf {\tiny 2}}
      \put(5,15){\textbf {\tiny \textit{l}}}
      \put(12,9){.}
      \put(10,10){.}
      \put(8,11){.}
      \put(0,0){\textbf {\tiny 0}}
    \end{picture}
  \end{minipage}
}

\def\sgraph{\unitlength.1em
  \begin{minipage}{25\unitlength}
    \begin{picture}(25,19)
      \put(5,2){\line(1,0){11}}
      \put(5,2){\line(0,1){11}}
      \put(5,2){\line(5,2){11}}
      \put(18,2){\textbf {\tiny 1}}
      \put(17,7){\textbf {\tiny 2}}
      \put(5,15){\textbf {\tiny \textit{s}}}
      \put(12,9){.}
      \put(10,10){.}
      \put(8,11){.}
      \put(0,0){\textbf {\tiny 0}}
    \end{picture}
  \end{minipage}
}

\begin{document}
\title{Floer homology and embedded bipartite graphs}
\author{Yuanyuan Bao}
\address{
Graduate School of Mathematical Sciences,
the University of Tokyo, Komaba 3-8-1, Meguro-ku, Tokyo 980-8577, Japan.
}

\email{bao@ms.u-tokyo.ac.jp}
\date{}
\begin{abstract}
We generalize the construction of the Heegaard Floer homology for a singular knot to that for a balanced bipartite graph. For a given graph, we provide a combinatorial description of the Euler characteristic of its Heegaard Floer homology by using the ``Kauffman states" on a graph diagram. 

\end{abstract}
\maketitle

\section{Introduction}

In \cite{MR2529302}, Ozsv{\'a}th, Stipsicz and Szab{\'o} studied the Heegaard Floer homology for a singular knot, which can be regarded as a four-valent spatial graph embedded in $S^{3}$. Their construction was latter used in \cite{MR2574747} by Ozsv{\'a}th and Szab{\'o} to provide an algebraic description of the Heegaard Floer homology for a knot. 

In this paper, we generalize the construction in \cite{MR2529302} to a balanced bipartite graph in $S^3$. In particular a trivalent graph without source or sink can be regarded as a bipartite graph, so our definition works for such a trivalent graph. For a balanced sutured manifold, Juh{\'a}sz \cite{MR2253454} defined sutured Floer homology. Alishahi and Eftekhary in \cite{MR3412088} provided a refinement of sutured Floer homology, where they assigned a variable to each suture of the manifold. For a balanced bipartite graph $G$ in $S^3$, its complement in $S^3$ is a sutured manifold whose sutures are the meridians of the edges of $G$. By applying Alishahi and Eftekhary's construction, we consider two chain complexes $\mathrm{CFG}^{\mathcal{V}}$ and $\mathrm{CFG}^{\mathcal{E}}$, by assigning a variable to each vertex (resp. edge) of $G$ in $\mathrm{CFG}^{\mathcal{V}}$ (resp. $\mathrm{CFG}^{\mathcal{E}}$). Harvey and O'Donnol \cite{harvey} recently constructed a combinatorial Heegaard Floer homology for a bipartite graph in $S^{3}$, where they defined the chain complex on a grid diagram of the given bipartite graph. If we regard a grid diagram as a special Heegaard diagram, the chain complex $\mathrm{CFG}^{\mathcal{V}}$ coincides with the chain complex defined there.

In the latter half of the paper, we study the Euler characteristic of $\widehat{\mathrm{HFG}}(G)$, which is the usual sutured Floer homology for the complement of the graph. The Euler characteristic of the sutured Floer homology has been well studied by Friedl, Juh{\'a}sz, and Rasmussen in \cite{MR2805998}. In this part, we provide a diagrammatical interpretation of the Euler characteristic of $\widehat{\mathrm{HFG}}(G)$.
In particular, we describe it combinatorially as a state sum over all `` Kauffman states" on the graph diagram. The description can be regarded as an extension of Kauffman's definition \cite{MR712133} for the Alexander polynomial of a link. When the graph is a $\theta_{n}$-curve, the construction is the Alexander polynomial studied in \cite{MR994083}.

In our subsequent papers \cite{alex, bao2}, we found that the state sum satisfies a series of relations which are analog of Murakami, Ohtsuki and Yamada's relations in \cite{MR1659228} for $U_q(\mathfrak{sl}(n))$-polynomial invariants when $n\geq 2$, and furthermore we proved that for a trivalent graph without source or sink the state sum 
is equivalent to the $\mathfrak{gl}(1\vert 1)$-Alexander polynomial of the trivalent graph that Viro defined in \cite{MR2255851}. Many interesting faces of the state sum are thus obtained.

Spatial graphs in the 3-sphere $S^{3}$ are widely used in the construction of quantum invariants for links and 3-manifolds, as shown in \cite{MR1036112}. For the fundamental representations of the quantum group $U_q(\mathfrak{sl}(n))$ ($n\geq 2$), Murakami, Ohtsuki and Yamada in \cite{MR1659228} defined an invariant for a trivalent planar graph, using which they provided a straightforward graphical construction of the associated quantum invariants of a link. The categorification of $\mathfrak{sl}(n)$-quantum invariants, which is now called $\mathfrak{sl}(n)$-link homology, has been actively studied by many authors. For $n=0$, such theory is the Heegaard Floer homology of a link constructed by Ozsv{\'a}th and Szab{\'o} \cite{MR2065507, MR2443092}, and independently Rasmussen \cite{Rasmussenthesis}, the construction of which in flavor is completely different from the $n>0$ cases. We hope the discussion in this paper may be useful in understanding the quantum topological side of Heegaard Floer theory, which is related to our next project.

\vspace{5pt}\noindent{\bf Acknowledgements.} This work was partially supported by Grant-in-Aid for Research Activity Start-up, and by Platform for Dynamic Approaches to Living 
System from the Ministry of Education, Culture, Sports, Science and 
Technology, Japan.

\section{Bipartite graphs embedded in $S^3$}
In this section, we define some basic terminologies and fix some notations for the graphs to be studied in this paper. 

\subsection{Balanced bipartite graphs}
\begin{defn}
\label{bipartite}
\rm
A graph $G$ with the vertex set $V$ and the edge set $E$ is called a \textit{bipartite graph} if $V$ is a disjoint union of two non-empty sets $V_{1}$ and $V_{2}$ so that each edge in $E$ is incident to both $V_{1}$ and $V_{2}$. If $\vert V_{1}\vert=\vert V_{2}\vert$, the graph $G$ is called \textit{balanced}.
\end{defn}

We assume that all the graphs in this paper have no isolated vertices or single-valent vertices.
The way of splitting $V=V_{1}\coprod V_{2}$ is not unique when $G$ is disconnected, while when $G$ is connected, it is very easy to see that the splitting is unique.

\begin{figure}
	\centering
		\includegraphics[width=0.4\textwidth]{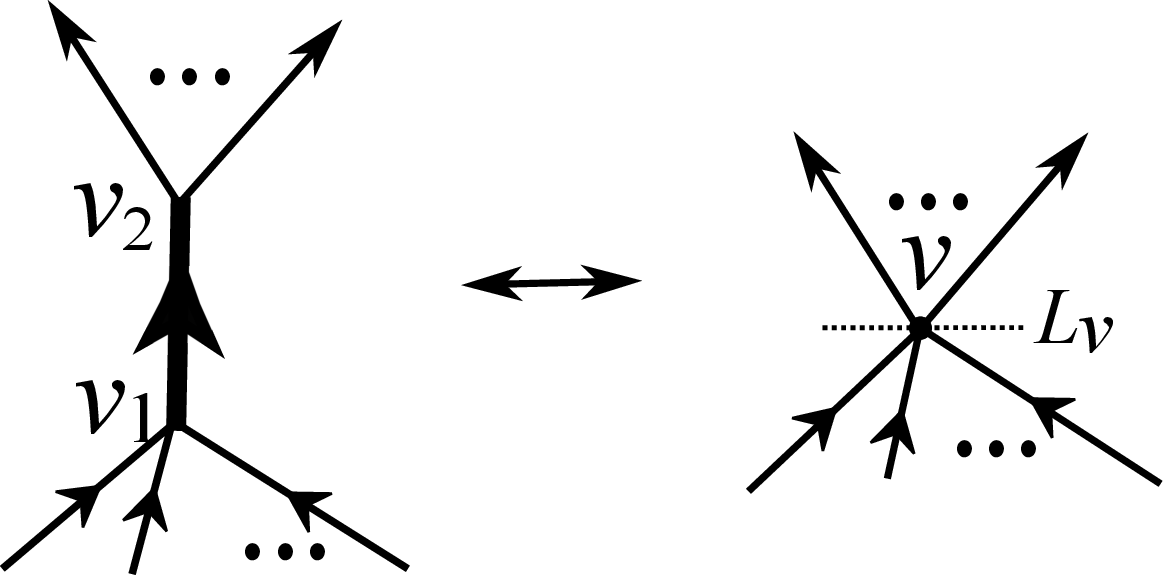}
	\caption{Thick edges characterize the orientation of $G$. After contracting the thick edges, we get a new graph where each new vertex $v$ is equipped with a disk $L_v$.}
	\label{fig:e32}
\end{figure}

Let $G$ be a balanced bipartite graph for which $n=\vert V_{1}\vert=\vert V_{2}\vert$. 
We consider an orientation of $G$ under which there are $n$ disjoint edges $\{e_{i}\}_{i=1}^{n}$ directing from $V_{1}$ to $V_{2}$ and the other edges direct from $V_{2}$ to $V_{1}$. Namely the set $\{e_{i}\}_{i=1}^{n}$ is a matching of $G$.
We call these $n$ distinguished edges the \emph{thick edges} of the oriented graph $G$. 

For a bipartite graph $G$ equipped with an orientation as above, we choose a diagram of $G$ in $S^2$ or $\mathbb{R}^2$
so that the thick edges are locally placed as in the left hand figure of Fig.~\ref{fig:e32}.  After contracting each of the thick edges to a vertex, we get a new graph, as on the right hand side of Fig.~\ref{fig:e32}. 

The graph on the right hand side satisfies the condition that at each vertex $v$, there is a small disk $L_v$ centered at $v$, which we call \emph{the dividing disk}, so that all the edges entering $v$ enter through one side of $L_v$ and all the edges leaving $v$ leave from the other side of $L_v$.  It is not hard to see that such a disk exists for an oriented graph if and only if the graph can be obtained from a balanced bipartite graph by contracting the thick edges as we described above. In light of this correspondence, we will use the diagram on the right hand side to represent an oriented balanced bipartite graph. 
It is named a transverse graph in \cite{harvey}. A balanced bipartite graph equipped with a balanced coloring is called an MOY graph in some papers, such as in H. Wu \cite{wu} and \cite{alex}.

\begin{rem}
Note that not every balanced bipartite graph allows the existence of thick edges. A simple example is a graph where there are two vertices in $V_1$ adjacent to one common vertex in $V_2$ but not to any other vertices. 
\end{rem}

\begin{figure}
	\centering
		\includegraphics[width=0.9\textwidth]{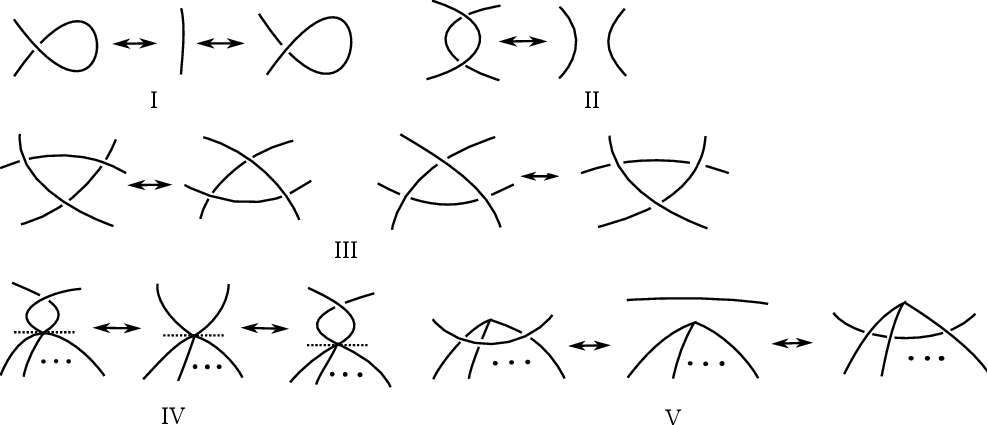}
	\caption{Reidemeister moves for balanced bipartite graphs. The twist in (IV) is not allowed to intersect the dividing disk $L_v$ at a vertex $v$.}
	\label{fig:e2}
\end{figure}

We have the following lemma, a proof of which can be found in \cite[Theorem 2.4]{harvey}.

\begin{lemma}
Two diagrams represent the same balanced bipartite graph if and only if they are connected by a sequence of Reidemeister moves in Fig. \ref{fig:e2}.
\end{lemma}

\subsection{Trivalent graphs}
Let $G$ be an oriented trivalent graph without source or sink embedded in $S^3$. The vertices of $G$ can be separated into two types as below.
\begin{align*}
\begin{tikzpicture}[baseline=-0.65ex, thick, scale=1.2]
\draw  (0, 0) [->-] to (0,0.5);
\draw [dotted]  (-0.5,0.5) to (0.5,0.5);
\draw   (0,0.5) [->]  to  (0.66,1);
\draw   (0,0.5) [->]  to  (-0.66,1);
\draw (0,-0.5) node {odd type};
\end{tikzpicture} \hspace{2cm}
\begin{tikzpicture}[baseline=-0.65ex, thick, scale=1.2]
\draw  (0, 0.5) [->] to (0,1);
\draw [dotted]  (-0.5,0.5) to (0.5,0.5);
\draw   (0,0.5) [-<-]  to  (0.66,0);
\draw   (0,0.5) [-<-]  to  (-0.66,0);
\draw (0,-0.5) node {even type};
\end{tikzpicture}
\end{align*}
The dividing disk automatically exists for a vertex of either odd or even type. A trivalent graph can be regarded as a bipartite graph, as we can see from the following lemma.

\begin{lemma}
Two trivalent graphs in $S^3$ are ambient isotopic to each other if and only if they are ambient isotopic as bipartite graphs.
\end{lemma}
\begin{proof}
It is easy to verify that for a trivalent graph, a Reidemeister move of type (IV) which does not respect the existence of the dividing disk at a vertex can be realized by a sequence of moves (I), (V) and moves of type (IV) that respect the dividing disk. 
\end{proof}


\section{Heegaard diagram for a bipartite graph}

\subsection{Heegaard diagram}
Consider an oriented balanced bipartite graph $G$ in $S^3$, which  is obtained from a usual bipartite graph by contracting a set of thick edges as we discussed before. Let $V$ denote the set of vertices of $G$ and $E$ denote the set of edges. Suppose $\vert V \vert =n$ and $\vert E \vert =m$. We define the Heegaard diagram for $G$ as follows.
\begin{defn}
\label{hd}
\rm
A quintet $(\Sigma, \boldsymbol{\alpha}, \boldsymbol{\beta}, \boldsymbol{w}, \boldsymbol{z})$ is called a Heegaard diagram for $G$ if it satisfies the following conditions.
\begin{enumerate}
\item $(\Sigma, \boldsymbol{\alpha}, \boldsymbol{\beta}, \boldsymbol{w})$ is an $n$-pointed Heegaard diagram for $S^{3}$, and $\boldsymbol{z}$ is a set of $m$ points in $\Sigma\setminus (\boldsymbol{\alpha} \cup \boldsymbol{\beta} \cup \boldsymbol{w})$. 

\item For each vertex $v\in V$ whose indegree is $l$ (resp. outdegree is $s$), there exists a smooth embedding $\varphi_{v}: (\,\lgraph \,, \{0\}, \lbrace 1, 2, \cdots, l\rbrace)\hookrightarrow (\Sigma \backslash \boldsymbol{\alpha}, \boldsymbol{w}, \boldsymbol{z})$ (resp. $\psi_{v}:(\,\sgraph \,,\{0\}, \lbrace 1, 2, \cdots, s\rbrace)\hookrightarrow (\Sigma \backslash \boldsymbol{\beta}, \boldsymbol{w}, \boldsymbol{z})$) so that
the images of $\varphi_{v}$ (resp. $\psi_{v}$) are pairwisely disjoint and $\bigcup_{v\in V}(\mathrm{Im}(\varphi_{v})\cup \mathrm{Im}(\psi_{v}))$ recovers $G$, where we push the interior of $\mathrm{Im}(\varphi_{v})$ (resp. $\mathrm{Im}(\psi_{v})$) slightly into ${U_{\boldsymbol{\alpha}}}$ (resp. ${U_{\boldsymbol{\beta}}}$). Here ${U_{\boldsymbol{\alpha}}}$ (resp. ${U_{\boldsymbol{\beta}}}$) is obtained from $\Sigma$ by attaching $2$-handles along $\boldsymbol{\alpha}$-curves (resp. $\boldsymbol{\beta}$-curves).
\end{enumerate}
\end{defn}

From the definition it is easy to see that each base point in $\boldsymbol{w}$ corresponds to a vertex of $G$, and each base point in $\boldsymbol{z}$ corresponds to an edge of $G$. In light of this correspondence, we let $\boldsymbol{w}=\{w_v\}_{v\in V}$ and $\boldsymbol{z}=\lbrace z_{e}\rbrace_{e\in E}$.

\begin{figure}[h]
	\centering
		\includegraphics[width=0.8\textwidth]{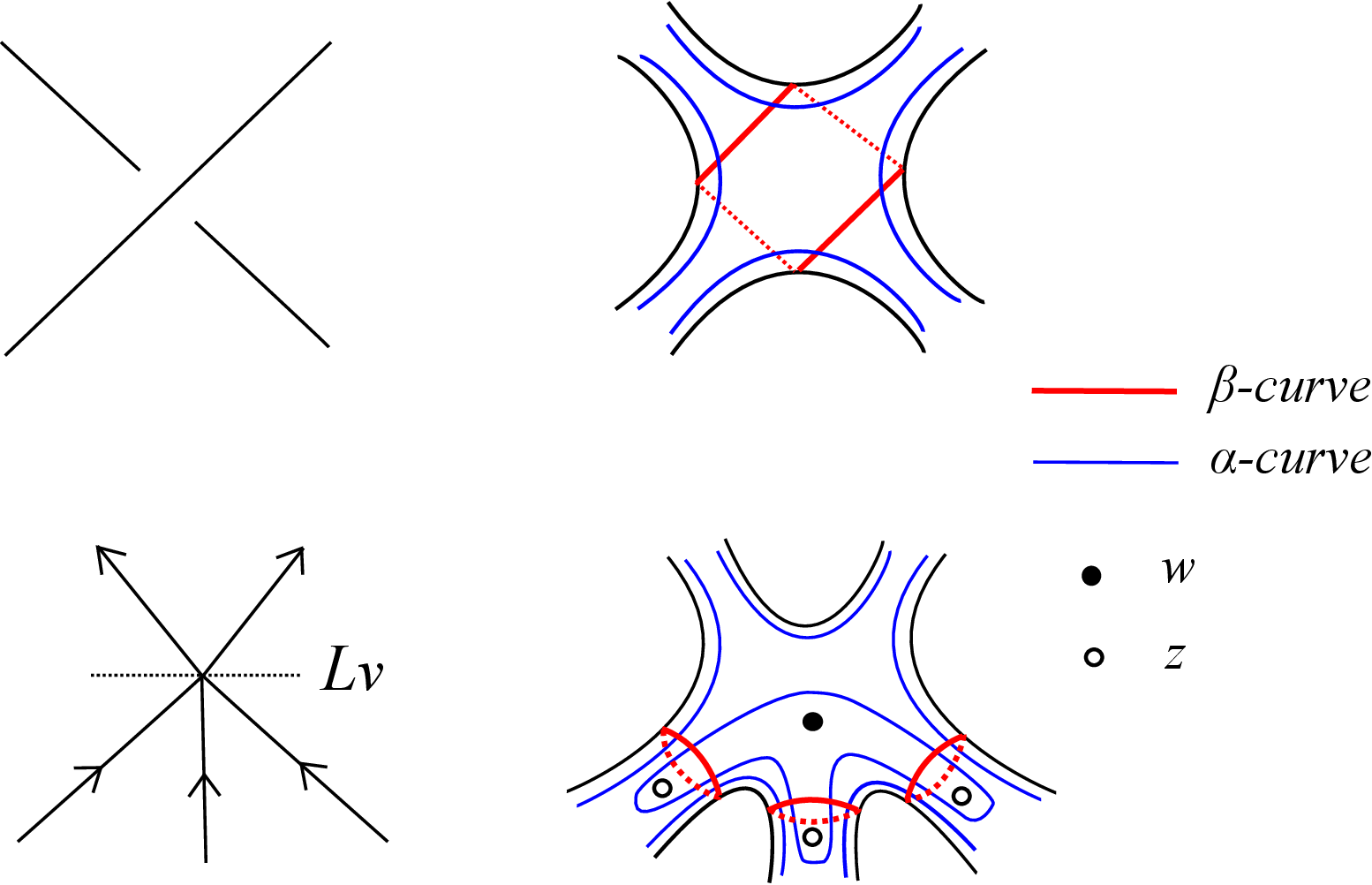}
		\caption{The Heegaard diagram associated with a graph diagram.}
	\label{fig:f1}
\end{figure}

\begin{ex}
\label{hdex}
\rm
We extend the ideas in \cite{MR1988285} and \cite{MR2529302} to provide a Heegaard diagram for a given balanced bipartite graph from its graph diagram in $S^{2}$. Consider a graph diagram $D\subset S^2$ for a given graph $G \subset S^{3}$. We assume that $D$ is connected as a projection.
\begin{enumerate}
\item Regard $D$ as a $1$-complex in $S^{3}$ and take a tubular neighbourhood of it in $S^{3}$. It is a handlebody and its boundary is the Heegaard surface $\Sigma$.
\item The diagram $D$ divides $S^{2}$ into several regions. For each region, introduce an $\boldsymbol{\alpha}$-curve on $\Sigma$ which encloses the region.
\item For each crossing of $D$, introduce a $\boldsymbol{\beta}$-curve following the rule in Figure~\ref{fig:f1}. 
\item Place the base point $w_v$ on each vertex $v\in V$.
\item Suppose a vertex $v$ has indegree $l$. introduce $l$ $\boldsymbol{\beta}$-curves which are meridians of the edges pointing to $v$ and $l$ base points of type $\boldsymbol{z}$ on the edges pointing to $v$. Introduce an $\boldsymbol{\alpha}$-curve $\alpha_v$ which bounds a disk around $w_v$ and encloses $w_v$ and all the base points of type $\boldsymbol{z}$ on the edges pointing to $v$. 

\item Remove one $\alpha_v$ and one $\boldsymbol{\alpha}$-curve created in Step (ii).

\end{enumerate}

It is easy to verify that the construction above gives a Heegaard diagram for $G$. 
\end{ex}

\subsection{Admissibility}
Let $(\Sigma, \boldsymbol{\alpha}, \boldsymbol{\beta}, \boldsymbol{w}, \boldsymbol{z})$ be a Heegaard diagram for a graph $G$ whose number of vertices is $n$. Let $$\mathbb{T}_{\boldsymbol{\alpha}}=\alpha_{1}\times \alpha_{1}\times \cdots \times \alpha_{n}  \text{ and }\mathbb{T}_{\boldsymbol{\beta}}=\beta_{1}\times \beta_{1}\times \cdots \times \beta_{n}$$ be the tori in the symmetric product $\mathrm{Sym}^{n}(\Sigma)$. Given $x, y\in \mathbb{T}_{\boldsymbol{\alpha}}\cap \mathbb{T}_{\boldsymbol{\beta}}$, let $\pi_{2}(x, y)$ be the set of relative homology classes of Whitney disks from $x$ to $y$ with boundary in $\mathbb{T}_{\boldsymbol{\alpha}}$ and $\mathbb{T}_{\boldsymbol{\beta}}$. For $\phi \in \pi_{2}(x, y)$, let $\mu(\phi)$ be its Maslov index and $\mathcal{\widehat{M}}(\phi)$ be the moduli space of pseudo-holomorphic disks in the class $\phi$ modulo $\mathbb{R}$.

Let $D_{1}, D_{2}, \cdots, D_{h}$ denote the closures of the components of $\Sigma\backslash (\boldsymbol{\alpha}\cup \boldsymbol{\beta})$. A \textit{domain} is a 2-chain on $\Sigma$ of the form $D=\sum_{i=1}^{h}a_{i}D_{i}$, where $a_{i}\in \mathbb{Z}$ is called the local multiplicity of $D$ at $D_{i}$. For a point $p$ in the interior of $D_i$, let $n_p(D)$ denote the local multiplicity of $D$ at the point $p$, which equals $a_i$. A domain $D$ is a \textit{positive domain} if $a_{i}\geq 0$ for $1\leq i\leq h$. 
A domain $P=\sum_{i=1}^{h}a_{i}D_{i}$ is called a \textit{periodic domain} if $\partial P$ is a $\mathbb{Z}$-linear combination of $\boldsymbol{\alpha}$-curves and $\boldsymbol{\beta}$-curves and $P\cap \boldsymbol{w}=P\cap \boldsymbol{z}=\emptyset$.

The Heegaard Floer complex is defined on a Heegaard diagram. In order for the differential to be well-defined, we need the following technical condition on the Heegaard diagram. 

\begin{defn}
\rm
A Heegaard diagram $(\Sigma, \boldsymbol{\alpha}, \boldsymbol{\beta}, \boldsymbol{w}, \boldsymbol{z})$ of a graph is said to be \textit{admissible} if every non-trivial periodic domain has both positive and negative local multiplicities at the Heegaard surface. 
\end{defn}

\begin{prop}
\label{graphadmi}
Suppose $G$ is a connected graph in $S^{3}$. Then the Heegaard diagram constructed in Example \ref{hdex} is admissible. 
\end{prop}
\begin{proof}
Let $A_v$ (resp. $B_v$) be the component of $\Sigma \backslash \boldsymbol{\alpha}$ (resp. $\Sigma \backslash \boldsymbol{\beta}$)
that contains the base point $w_v$ of $G$. Then any periodic domain has the form $$P=\sum_{v\in V} a_{v}A_{v}+b_{v}B_{v},$$ where $a_{v}, b_{v}\in \mathbb{Z}$. Since $n_{w_{v}}(P)=0$, we have $a_{v}+b_{v}=0$ for any $v\in V$. Let $e_{uv}$ be an edge of $G$ that directs from the vertex $u$ to $v$. Then we have $a_v+b_u=0$ since $n_{z_{e_{uv}}}(P)=0$, which implies that $a_v=a_u$ and $b_v=b_u$.

Since $G$ is connected, any two vertices of $G$ is connected by a sequence of edges. Therefore we have $a_v=a_u$ and $b_v=b_u$ for any vertices $u, v$ of $G$. As a result we see that $P$ is the trivial domain.  
\end{proof}

\begin{rem}
For a disconnected graph, the Heegaard diagram constructed in Example \ref{hdex} is not admissible.
\end{rem}

\section{Heegaard Floer Complex for a bipartite graph}
It is easy to see that the complement of a balanced bipartite graph is a balanced sutured manifold where the sutures are given by the meridians of the edges and the meridian circles around vertices. The Heegaard Floer homology for a sutured manifold and its refinement have been constructed in \cite{MR2253454} and \cite{MR3412088}. In this section, we interpret how their theories can be applied to the case of bipartite graphs to extract interesting invariants for a graph.  

\subsection{Alishahi-Eftekhary's refinement}
For a balanced sutured manifold, Juh{\'a}sz \cite{MR2253454} defined the sutured Floer homology for it, which is defined on a Heegaard diagram where each suture corresponds to a base point. Alishahi and Eftekhary \cite{MR3412088} extended Juh{\'a}sz's definition and provided a minus version of sutured Floer homology, where they assigned a variable to each of the sutures. 

We briefly review their definition. For details, please refer to the original papers. Let $(X, \tau)$ be a balanced sutured manifold and $\tau=\{\gamma_1, \gamma_2, \cdots, \gamma_k\}$ be the set of sutures. The set $\tau$ divides $\partial X$ into two parts. Namely $\partial X-\tau=\mathfrak{R}^{+}(\tau)\cup \mathfrak{R}^{-}(\tau)$, where $\mathfrak{R}^{+}(\tau)$ and $\mathfrak{R}^{-}(\tau)$ are called the positive and the negative part respectively. Let $\lambda_j$ be the variable associated with $\gamma_j$ for $1\leq j \leq k$.

Let $\{R_1^{+}, R_2^{+}, \cdots, R_p^{+} \}$ be the connected components of $\mathfrak{R}^{+}(\tau)$. Let $$\lambda_i^{+}=\prod_{\gamma_j \subset \partial R_i^{+}} \lambda_j \text{ for $1\leq i \leq p$, and } \lambda^{+}(\tau)=\sum_{i=1}^{p}  \lambda_i^{+} .$$
Similarly one can define $\lambda^{-}(\tau)$. Let  
$$\mathbb{A}_{\tau}=\frac{\langle \lambda_1, \lambda_2, \cdots, \lambda_k \rangle_{\mathbb{F}}}{\langle \lambda^{+}(\tau)-\lambda^{-}(\tau)\rangle_{\mathbb{F}} + V^{+} + V^{-}},$$ where $\langle \lambda_1, \lambda_2, \cdots, \lambda_k \rangle_{\mathbb{F}}$ is the polynomial ring over $\mathbb{F}:=\mathbb{Z}/ 2\mathbb{Z}$ generated by $\lambda_j$ for $1\leq j \leq k$ and $V^{+}$ and $V^{-}$ are two subalgebras of $\langle \lambda_1, \lambda_2, \cdots, \lambda_k \rangle_{\mathbb{F}}$ defined from the components of $\mathfrak{R}^{+}(\tau)$ and $\mathfrak{R}^{-}(\tau)$ with positive genera.

Given a Heegaard diagram for $(X, \tau)$. The chain complex $\mathrm{CF}(X, \tau)$ in \cite{MR3412088} is generated as a free $\mathbb{A}_{\tau}$-module by $\mathbb{T}_{\boldsymbol{\alpha}}\cap \mathbb{T}_{\boldsymbol{\beta}}$, and the differential is defined as follows.
\begin{eqnarray*}
\begin{split}
\partial(x)
=\sum_{y\in \mathbb{T}_{\boldsymbol{\alpha}}\cap \mathbb{T}_{\boldsymbol{\beta}}} \sum_{\{\phi \in \pi_{2}(x, y)| \mu (\phi)=1\}} \sharp \mathcal{\widehat{M}}(\phi)\cdot \prod_{j=1}^{k} \lambda_{j}^{n_{z_j}(\phi)}y
\end{split}
\end{eqnarray*}
for any $x\in \mathbb{T}_{\boldsymbol{\alpha}}\cap\mathbb{T}_{\boldsymbol{\beta}}$, where $z_j$ is the base point associated with $\gamma_j$ for $1\leq j \leq k$.

A relative $H_1(X; \mathbb{Z})$-filtration on $\mathrm{CF}(X, \tau)$ is determined by the following rules.
\begin{eqnarray*}
A(\lambda_j)&=&[\gamma_j]  \text{ for $1\leq j \leq k$, and}\\
A(x)-A(y)&=& A(\prod_{j=1}^{k} \lambda_{j}^{n_{z_j}(\phi)})=\sum_{j=1}^{k} n_{z_j}(\phi) [\gamma_j],
\end{eqnarray*}
for any $x, y\in \mathbb{T}_{\boldsymbol{\alpha}}\cap\mathbb{T}_{\boldsymbol{\beta}}$ and $\phi \in \pi_2(x, y)$.

\begin{theo}[\cite{MR3412088}]
\label{bigtheo}
The filtered chain homotopy type of the $\mathbb{A}_{\tau}$-chain complex $\mathrm{CF}(X, \tau)$ with $H_1(X; \mathbb{Z})$-filtration is a topological invariant of the balanced sutured manifold $(X, \tau)$.
\end{theo}

\subsection{The chain complex $\mathrm{CFG}^{\mathcal{V}}$ for a graph: assign a variable to each vertex}
\begin{figure}
	\centering
		\includegraphics[width=0.5\textwidth]{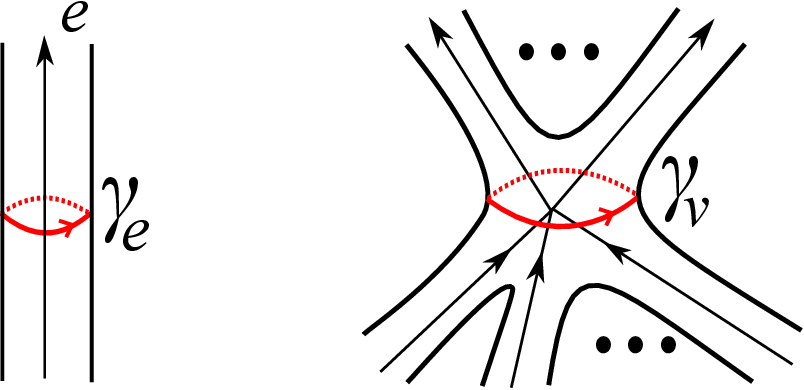}
		\caption{The meridian circle for an edge and for a vertex.}
	\label{fig:e11}
\end{figure}

For a balanced bipartite graph $G$, let $X$ be its complement in $S^3$. Then $(X, \tau=\{\gamma_{e}\}_{e\in E} \cup \{\gamma_{v}\}_{v\in V})$ becomes a balanced sutured manifold where $\gamma_{e}$ denotes the oriented meridian of the edge $e$ and $\gamma_{v}$ is the oriented meridian circle around the vertex $v$, as shown in Fig. \ref{fig:e11}. In this case we see that $V^{+}$ and $V^{-}$ in Alishahi-Eftekhary's construction vanish for the reason that there exists no positive genus component in $\mathfrak{R}^{+}(\tau)$ or $\mathfrak{R}^{-}(\tau)$.

We assign a variable to each vertex of $G$ and zero to each edge of $G$. The algebra $\mathbb{A}_{\tau}$ becomes a free commutative algebra since both 
$\lambda^{+}(\tau)$ and $\lambda^{-}(\tau)$ vanish. We obtain a version of Floer homology for a graph. More precisely, let $U_v$ be the variable assigned to the vertex $v$ of $G$, and let $\mathbb{F}[\{U_{v}\}_{v\in V}]$ be the polynomial ring over $\mathbb{F}$ generated by $\{U_{v}\}_{v\in V}$. Suppose $(\Sigma, \boldsymbol{\alpha}, \boldsymbol{\beta}, \boldsymbol{w}, \boldsymbol{z})$ is an admissible Heegaard diagram for $G$. We define the chain complex $(\mathrm{CFG}^{\mathcal{V}}(\Sigma, \boldsymbol{\alpha}, \boldsymbol{\beta}, \boldsymbol{w}, \boldsymbol{z}), \partial^{\mathcal{V}})$ as follows. It is a free $\mathbb{F}[\{U_{v}\}_{v\in V}]$-module generated by $\mathbb{T}_{\boldsymbol{\alpha}}\cap \mathbb{T}_{\boldsymbol{\beta}}$ and the differential is
\begin{eqnarray*}
\begin{split}
\partial^{\mathcal{V}}(x)
=\sum_{y\in \mathbb{T}_{\boldsymbol{\alpha}}\cap \mathbb{T}_{\boldsymbol{\beta}}} \sum_{\{\phi \in \pi_{2}(x, y)| \mu (\phi)=1, n_{\boldsymbol{z}}(\phi)=\{0\}\}} \sharp \mathcal{\widehat{M}}(\phi)\cdot \prod_{v\in V} U_{v}^{n_{w_{v}}(\phi)}y
\end{split}
\end{eqnarray*}
for any $x\in \mathbb{T}_{\boldsymbol{\alpha}}\cap\mathbb{T}_{\boldsymbol{\beta}}$, where $n_{\boldsymbol{z}}(\phi)=\{n_{z}(\phi)\}_{z\in \boldsymbol{z}}$. The differential here only counts those pseudo-holomorphic disks avoiding the base points of type $\boldsymbol{z}$.

The relative Maslov grading of the complex is defined by the following relations
\begin{eqnarray*}
M(x)-M(y)&=&\mu(\phi)-2\sum_{v\in V}n_{w_{v}}(\phi) \text{ and }\\ 
M(U_{v})&=&-2,
\end{eqnarray*}
for any $x, y\in \mathbb{T}_{\boldsymbol{\alpha}}\cap\mathbb{T}_{\boldsymbol{\beta}}$, where $\phi$ is any Whitney disk from $x$ to $y$.

The first homology group $H_{1}(X; \mathbb{Z})$ is generated by $[\gamma_e]$'s with the relation $$\prod_{\text{$e$ pointing to $v$}}[\gamma_e]=\prod_{\text{$e$ pointing out of $v$}}[\gamma_e],$$ for any vertex $v$ of $G$. It is easy to check that $[\gamma_{v}]$ equals the value above.
The chain complex above has a relative $H_{1}(X; \mathbb{Z})$-grading which we call the Alexander grading. It is defined by the following relations
\begin{eqnarray*}
A(x)-A(y)&=&\prod_{v\in V}[\gamma_v]^{n_{w_v}(\phi)}\div \prod_{e\in E}[\gamma_e]^{n_{z_e}(\phi)}, \text{ and }\\
A(U_{v})&=&[\gamma_{v}],
\end{eqnarray*}
for any $x, y\in \mathbb{T}_{\boldsymbol{\alpha}}\cap\mathbb{T}_{\boldsymbol{\beta}}$, where $\phi$ is any Whitney disk from $x$ to $y$. Unlike Alishahi and Eftekhary's convention, here we take the multiplication as the group operation in $H_{1}(X; \mathbb{Z})$ for the sake of convenience of discussion in Section 5.

\begin{lemma}
The differential $\partial^{\mathcal{V}}$ decreases the Maslov grading by one and preserves the Alexander grading. 
\end{lemma}
\begin{proof}
Follow from the definition of $\partial^{\mathcal{V}}$.
\end{proof}

As a special case of Theorem \ref{bigtheo}, we have the following result.
\begin{prop}
The homology of $\mathrm{CFG}^{\mathcal{V}}(\Sigma, \boldsymbol{\alpha}, \boldsymbol{\beta}, \boldsymbol{w}, \boldsymbol{z})$, which is denoted by $$\mathrm{HFG}^{\mathcal{V}}(G)=\bigoplus_{j\in \mathbb{Z}, \alpha\in H_{1}(X; \mathbb{Z})}\mathrm{HFG}^{\mathcal{V}}_j(G; \alpha),$$ is a topological invariant of $G$, where $j$ denotes the Maslov grading and $\alpha$ is the Alexander grading.
\end{prop}

\subsection{The chain complex $\mathrm{CFG}^{\mathcal{E}}$ for a graph: assign a variable to each edge}
If we assign a variable to each edge of the graph $G$ and zero to each vertex, we obtain another version of Floer homology for the graph. For a singular knot, which can be regarded as a balanced bipartite graph, Ozsv{\'a}th and Szab{\'o} \cite{MR2529302, MR2574747} studied the Heegaard Floer homology for it, where they assigned a variable to each edge of the singular knot. The following chain complex is an extension of their construction.

Suppose $(\Sigma, \boldsymbol{\alpha}, \boldsymbol{\beta}, \boldsymbol{w}, \boldsymbol{z})$ is an admissible Heegaard diagram for $G$. We define the chain complex $(\mathrm{CFG}^{\mathcal{E}}(\Sigma, \boldsymbol{\alpha}, \boldsymbol{\beta}, \boldsymbol{w}, \boldsymbol{z}), \partial^{\mathcal{E}})$ as below. It is a free $\mathbb{F}[\{U_{e}\}_{e\in E}]$-module generated by $\mathbb{T}_{\boldsymbol{\alpha}}\cap \mathbb{T}_{\boldsymbol{\beta}}$ and the differential is
\begin{eqnarray*}
\begin{split}
\partial^{\mathcal{E}}(x)
=\sum_{y\in \mathbb{T}_{\boldsymbol{\alpha}}\cap \mathbb{T}_{\boldsymbol{\beta}}} \sum_{\{\phi \in \pi_{2}(x, y)| \mu (\phi)=1, n_{\boldsymbol{w}}(\phi)=\{0\}\}} \sharp \mathcal{\widehat{M}}(\phi)\cdot \prod_{e\in E} U_{e}^{n_{z_{e}}(\phi)}y,
\end{split}
\end{eqnarray*}
for any $x\in \mathbb{T}_{\boldsymbol{\alpha}}\cap\mathbb{T}_{\boldsymbol{\beta}}$. The differential here only counts those pseudo-holomorphic disks avoiding the base points in $\boldsymbol{w}$.

The relative Maslov grading of the complex is defined by the following relations
\begin{eqnarray*}
M(x)-M(y)&=&\mu(\phi)-2\sum_{v\in V}n_{w_{v}}(\phi) \text{ and }\\ 
M(U_{e})&=&0,
\end{eqnarray*}
for any $x, y\in \mathbb{T}_{\boldsymbol{\alpha}}\cap\mathbb{T}_{\boldsymbol{\beta}}$, where $\phi$ is any Whitney disk from $x$ to $y$. 

The chain complex above also has a relative $H_{1}(X; \mathbb{Z})$-grading which we call the Alexander grading. It is defined by the following relations
\begin{eqnarray*}
A(x)-A(y)&=&\prod_{v\in V}[\gamma_v]^{n_{w_v}(\phi)}\div \prod_{e\in E}[\gamma_e]^{n_{z_e}(\phi)}, \text{ and }\\
A(U_{e})&=&[\gamma_{e}]^{-1}
\end{eqnarray*}
for any $x, y\in \mathbb{T}_{\boldsymbol{\alpha}}\cap\mathbb{T}_{\boldsymbol{\beta}}$, where $\phi$ is any Whitney disk from $x$ to $y$. 

\begin{lemma}
The differential $\partial^{\mathcal{E}}$ decreases the Maslov grading by one and preserves the Alexander grading.
\end{lemma}
\begin{proof}
The variable $U_e$ does not affect the relative Maslov grading. The lemma follows from the definition of $\partial^{\mathcal{E}}$.
\end{proof}

As a special case of Theorem \ref{bigtheo}, we have the following result.

\begin{prop}
The homology of $\mathrm{CFG}^{\mathcal{E}}(\Sigma, \boldsymbol{\alpha}, \boldsymbol{\beta}, \boldsymbol{w}, \boldsymbol{z})$, which is denoted by $$\mathrm{HFG}^{\mathcal{E}}(G)=\bigoplus_{j\in \mathbb{Z}, \alpha\in H_{1}(X; \mathbb{Z})}\mathrm{HFG}^{\mathcal{E}}_j(G; \alpha),$$ is a topological invariant of $G$, where $j$ is the Maslov grading and $\alpha$ is the Alexander grading.
\end{prop}

\subsection{Hat version}
Let $(\widehat{\mathrm{CFG}}(\Sigma, \boldsymbol{\alpha}, \boldsymbol{\beta}, \boldsymbol{w}, \boldsymbol{z}), \widehat{\partial})$ be the chain complex obtained from $\mathrm{CFG}^{\mathcal{V}}(\Sigma, \boldsymbol{\alpha}, \boldsymbol{\beta}, \boldsymbol{w}, \boldsymbol{z})$ or $\mathrm{CFG}^{\mathcal{E}}(\Sigma, \boldsymbol{\alpha}, \boldsymbol{\beta}, \boldsymbol{w}, \boldsymbol{z})$ by setting all the variables to zero. It is the usual sutured Floer homology of the complement of $G$. It is easy to see that the Maslov gradings in Sections 4.2 and 4.3 induce the same grading in $\widehat{\mathrm{CFG}}(\Sigma, \boldsymbol{\alpha}, \boldsymbol{\beta}, \boldsymbol{w}, \boldsymbol{z})$, and the same holds true for the Alexander gradings. The homology group $\widehat{\mathrm{HFG}}(G)$ is a topological invariant of $G$.  

\subsection{Basic symmetries}
\begin{prop}

Given an oriented bipartite graph $G$, let $-G$ be the same graph with the reverse orientation. Then we have $$\mathrm{HFG}^{\omega}_{j}(G; \alpha) =\mathrm{HFG}^{\omega}_{j}(-G; \alpha),$$
for $\omega=\mathcal{V}, \mathcal{E}$.
\end{prop}

\begin{proof}
If $(\Sigma, \boldsymbol{\alpha}, \boldsymbol{\beta},\boldsymbol{w}, \boldsymbol{z})$ is a Heegaard diagram for $G$, then $(-\Sigma, \boldsymbol{\beta}, \boldsymbol{\alpha}, \boldsymbol{w}, \boldsymbol{z})$ is a Heegaard diagram for $-G$. For any pesudo-holomorphic disk $\phi \in \pi_{2}(x, y)$ for $x, y \in \mathbb{T}_{\boldsymbol{\alpha}}\cap \mathbb{T}_{\boldsymbol{\beta}}$ in $\Sigma$, we see that $-\phi$ is a pesudo-holomorphic disk in $-\Sigma$ connecting $x$ to $y$ for $x, y \in \mathbb{T}_{\boldsymbol{\beta}}\cap \mathbb{T}_{\boldsymbol{\alpha}}$. Therefore, there exists a natural chain isomorphism between $\mathrm{CFG}^{\omega}(\Sigma, \boldsymbol{\alpha}, \boldsymbol{\beta}, \boldsymbol{w}, \boldsymbol{z})$ and $\mathrm{CFG}^{\omega}(-\Sigma, \boldsymbol{\beta}, \boldsymbol{\alpha}, \boldsymbol{w}, \boldsymbol{z})$. The Alexander grading keeps invariant since we take the identification $H_1(S^3\backslash G; \mathbb{Z})=H_1(S^3\backslash (-G); \mathbb{Z})$ by sending $\gamma_e$ to the same curve with the reverse orientation.
\end{proof}

\begin{prop}
Given a graph $G$, let $G^{*}$ be its mirror image.  Then $$\widehat{\mathrm{HFG}}_{j}(G^{*}; \alpha) = \widehat{\mathrm{HFG}}_{-j}(G; -\alpha).$$. 

\end{prop}

\begin{proof}
If $(\Sigma, \boldsymbol{\alpha}, \boldsymbol{\beta},\boldsymbol{w}, \boldsymbol{z})$ is a Heegaard diagram for $G$, then $(-\Sigma, \boldsymbol{\alpha}, \boldsymbol{\beta}, \boldsymbol{w}, \boldsymbol{z})$ is a Heegaard diagram for $G^{*}$. For any pesudo-holomorphic disk $\phi \in \pi_{2}(x, y)$ for $x, y \in \mathbb{T}_{\boldsymbol{\alpha}}\cap \mathbb{T}_{\boldsymbol{\beta}}$ in $\Sigma$, it is easy to see that $-\phi$ is a pesudo-holomorphic disk in $-\Sigma$ connecting $y$ to $x$ for $x, y \in \mathbb{T}_{\boldsymbol{\alpha}}\cap \mathbb{T}_{\boldsymbol{\beta}}$. Therefore $\widehat{\mathrm{CFG}}(-\Sigma, \boldsymbol{\alpha}, \boldsymbol{\beta}, \boldsymbol{w}, \boldsymbol{z})$ is the dual complex of $\widehat{\mathrm{CFG}}(\Sigma, \boldsymbol{\alpha}, \boldsymbol{\beta}, \boldsymbol{w}, \boldsymbol{z})$. 

Suppose $x^*$ is the dual of $x$ in $\widehat{\mathrm{CFG}}(-\Sigma, \boldsymbol{\alpha}, \boldsymbol{\beta}, \boldsymbol{w}, \boldsymbol{z})$ for $x \in \mathbb{T}_{\boldsymbol{\alpha}}\cap \mathbb{T}_{\boldsymbol{\beta}}$. Note that the relative Maslov grading and Alexander grading in $\widehat{\mathrm{CFG}}(-\Sigma, \boldsymbol{\alpha}, \boldsymbol{\beta}, \boldsymbol{w}, \boldsymbol{z})$ are defined by the relations $M^{*}(x^*)-M^{*}(y^*)=M(y)-M(x)$ and $A^{*}(x^*)-A^{*}(y^*)=A(y)-A(x)$ for $x, y \in \mathbb{T}_{\boldsymbol{\alpha}}\cap \mathbb{T}_{\boldsymbol{\beta}}$. Therefore we have $\widehat{\mathrm{HFG}}_{j}(G^{*}; \alpha) = \widehat{\mathrm{HFG}}^{-j}(G; -\alpha)$. On the other hand, since we are working over the field $\mathbb{F}=\mathbb{Z}/2\mathbb{Z}$, we have $\widehat{\mathrm{HFG}}^{-j}(G; -\alpha)=\widehat{\mathrm{HFG}}_{-j}(G; -\alpha)$.

\end{proof}

\section{Euler characteristic}
In this section, we study the Euler characteristic of $\widehat{\mathrm{HFG}}(G)$ for a graph $G$ in $S^3$. 
We first provide a definition of the Alexander invariant for a graph, by extending Litherland's definition for a $\theta_{n}$-graph in \cite{MR994083}. 
Then we discuss the calculation of the invariant based on Fox calculus, which shows that the invariant coincides with the Euler characteristic.   

For a sutured manifold, Friedl, Juh{\'a}sz and Rasmussen in \cite{MR2805998} defined a torsion invariant and showed that it is the Euler characteristic of the sutured Floer homology of the underlying manifold. The discussion in this section can be regarded as a combinatorial interpretation of their result in the case of bipartite graphs. In particular, we provide a state sum formula for the Alexander invariant, which is a diagrammatical formula constructed on a graph diagram of $G$. 

\subsection{Euler characteristic}

Let $(\Sigma, \boldsymbol{\alpha}, \boldsymbol{\beta}, \boldsymbol{w}, \boldsymbol{z})$ be a Heegaard diagram for a given bipartite graph $G$. 
\begin{defn}
\rm
The \textit{Euler characteristic} for $\widehat{\mathrm{CFG}}(\Sigma, \boldsymbol{\alpha}, \boldsymbol{\beta}, \boldsymbol{w}, \boldsymbol{z})$ is 
$$\chi(\widehat{\mathrm{CFG}}(\Sigma, \boldsymbol{\alpha}, \boldsymbol{\beta}, \boldsymbol{w}, \boldsymbol{z})):=\sum_{x\in \mathbb{T}_{\boldsymbol{\alpha}}\cap \mathbb{T}_{\boldsymbol{\beta}}} (-1)^{M(x)}\cdot A(x), $$ which is an element in $\mathbb{Z} H_{1}(X; \mathbb{Z})$, where $X$ is the complement of $G$. 
\end{defn}

Since both $M(x)$ and $A(x)$ are relative gradings, the Euler characteristic is only well-defined modulo $\pm H_{1}(X; \mathbb{Z})$. The value $(-1)^{M(x)}$ gives a relative $\mathbb{Z}/2\mathbb{Z}$-grading on $\widehat{\mathrm{CFG}}(\Sigma, \boldsymbol{\alpha}, \boldsymbol{\beta}, \boldsymbol{w}, \boldsymbol{z})$, which can be calculated by considering the sign of each generator as below. The Heegaard surface $\Sigma$ has an orientation inherited from that of $S^3$. We choose orientations for the $\boldsymbol{\alpha}$-curves $\alpha_{1}, \alpha_{2}, \cdots, \alpha_{n}$ and the $\boldsymbol{\beta}$-curves $\beta_{1}, \beta_{2}, \cdots, \beta_{n}$. 
Identify $x\in \mathbb{T}_{\boldsymbol{\alpha}}\cap \mathbb{T}_{\boldsymbol{\beta}}$ with $\{x_{1}, x_{2}, \cdots, x_{n}\}$, where $x_{i}\in \alpha_{i}\cap \beta_{\sigma (i)}$ for some $\sigma$ in the symmetric group $S_{n}$. Denote the sign of the intersection point $x_{i}$ in $\Sigma$ by $\mathrm{sign}(x_{i})$. Then we consider $$\mathrm{sign}(x):=\mathrm{sign}(\sigma)\prod_{i=1}^{n}\mathrm{sign}(x_{i})\in \{1, -1\}.$$
The differential of $\widehat{\mathrm{CFG}}(\Sigma, \boldsymbol{\alpha}, \boldsymbol{\beta}, \boldsymbol{w}, \boldsymbol{z})$ changes $\mathrm{sign}(x)$, so we have $\mathrm{sign}(x)=(-1)^{M(x)}$ up to an overall sign change. As a result we see that 
\begin{equation}
\label{euler}
\chi(\widehat{\mathrm{CFG}}(\Sigma, \boldsymbol{\alpha}, \boldsymbol{\beta}, \boldsymbol{w}, \boldsymbol{z})) \doteq \sum_{x\in \mathbb{T}_{\boldsymbol{\alpha}}\cap \mathbb{T}_{\boldsymbol{\beta}}} \mathrm{sign}(\sigma)\prod_{i=1}^{n}\mathrm{sign}(x_{i})\cdot A(x),
\end{equation}
where $\doteq$ means the values which are connected are equal modulo $\pm H_{1}(X; \mathbb{Z})$.

\subsection{The Alexander invariant for a balanced bipartite graph}
\begin{figure}
	\centering
		\includegraphics[width=0.25\textwidth]{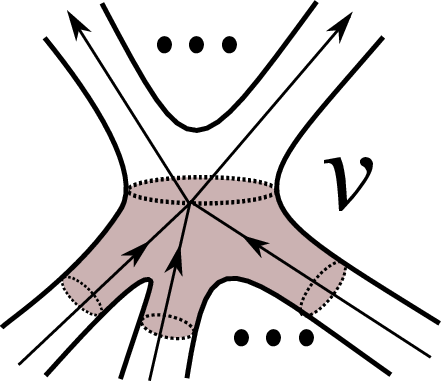}
		\caption{Shadowed surface is $\partial_{\mathrm{in}}(v)$.}
	\label{e10}
\end{figure}

Given an oriented balanced bipartite graph $G$ in $S^3$, let $X$ be its complement in $S^3$. Around each vertex $v$, consider a subsurface of $\partial X$ which is bounded by the meridians of the edges pointing to $v$ and the meridian circle around $v$, as the shadowed part in Figure~\ref{e10}. We call this subsurface $\partial_{\mathrm{in}}(v)$ and let $\partial_{\mathrm{in}}(X)$ be the disjoint union of $\partial_{\mathrm{in}}(v)$'s for all the vertices $v\in V$.

Let $\cal{G}=\mathbb{Z}H_{1}(X; \mathbb{Z})$, which is a commutative ring with unit. 
Let $\widetilde{X}$ be the universal abelian cover of $X$ defined by the Hurewicz map from $\pi_{1}(X, x_{0})$ to $H_{1}(X; \mathbb{Z})$, and $\widetilde{\partial}_{\mathrm{in}}(X)$ be the pre-image of $\partial_{\mathrm{in}}(X)$ in $\widetilde{X}$. The deck transformation group of $\widetilde{X}$ is isomorphic to $H_{1}(X; \mathbb{Z})$, so $H_{1}(\widetilde{X}, \widetilde{\partial}_{\mathrm{in}}(X); \mathbb{Z})$ has a $\cal{G}$-module structure. We call the $\cal{G}$-module $H_{1}(\widetilde{X}, \widetilde{\partial}_{\mathrm{in}}(X))$ the {\it Alexander module} of $G$, which is denoted by $A(G)$. Obviously, the isomorphism class of $A(G)$ is a topological invariant of $G$. When $G$ is a $\theta_n$-graph, this invariant has been studied by Litherland in \cite{MR994083}.

We show how to study $A(G)$ from a Heegaard diagram of $G$. Suppose $H_{G}=(\Sigma, \boldsymbol{\alpha}=\lbrace\alpha_{1}, \alpha_{2}, \cdots, \alpha_{n}\rbrace, \boldsymbol{\beta}=\lbrace\beta_{1}, \beta_{2}, \cdots, \beta_{n}\rbrace, \boldsymbol{w}, \boldsymbol{z})$ is a Heegaard diagram for the graph $G$. Then the Heegaard diagram $H_{G}$ provides a relative handle decomposition of $X$ built on $\partial_{\mathrm{in}}(X)$. We first attach 1-handles to $\partial_{\mathrm{in}}(X)\times I$ with belt circles $\alpha_{1}, \alpha_{2}, \cdots, \alpha_{n}$, and then attach 2-handles  with attaching circles $\beta_{1}, \beta_{2}, \cdots, \beta_{n}$. Then the relative chain complex $C_{1}(X, \partial_{\mathrm{in}}(X); \mathbb{Z})$ is freely generated by $\alpha_{1}^{*}, \alpha_{2}^{*}, \cdots, \alpha_{n}^{*}$, where $\alpha_{i}^{*}$ is the 1-handle with belt circle $\alpha_{i}$, and $C_{2}(X, \partial_{\mathrm{in}}(X); \mathbb{Z})$ is freely generated by $\beta_{1}^{*}, \beta_{2}^{*}, \cdots, \beta_{n}^{*}$, where $\beta_{i}^{*}$ is the 2-handle with attaching circle $\beta_{i}$.

The pre-images (lifts) of $\{\beta_{i}^{*}\}_{i=1}^{n}$ and those of $\{\alpha_{i}^{*}\}_{i=1}^{n}$ in $\widetilde{X}$ generate $C_{2}(\widetilde{X}, \widetilde{\partial}_{\mathrm{in}}(X); \mathbb{Z})$ and $C_{1}(\widetilde{X}, \widetilde{\partial}_{\mathrm{in}}(X); \mathbb{Z})$ as $\mathbb{Z}$-modules, respectively. Choose a lift $\widetilde{\beta_{i}^{*}}$ of $\beta_{i}^{*}$, and a lift $\widetilde{\alpha_{i}^{*}}$ of $\alpha_{i}^{*}$ in $\widetilde{X}$, where $i=1,2,\cdots, n$. Then $\lbrace \widetilde{\beta_{i}^{*}}\rbrace _{i=1}^{n}$ and $\lbrace \widetilde{\alpha_{i}^{*}}\rbrace _{i=1}^{n}$ generate $C_{2}(\widetilde{X}, \widetilde{\partial}_{\mathrm{in}}(X);\cal{G})$ and $C_{1}(\widetilde{X}, \widetilde{\partial}_{\mathrm{in}}(X); \cal{G})$ as free $\cal{G}$-modules, respectively. The relative chain complex for $(\widetilde{X}, \widetilde{\partial}_{\mathrm{in}}(X))$ defined in this way is zero except in dimension 1 and 2. We have
$$0  \to C_{2}(\widetilde{X}, \widetilde{\partial}_{\mathrm{in}}(X); \cal{G}) \xrightarrow{\widetilde{\partial}}  C_{1}(\widetilde{X}, \widetilde{\partial}_{\mathrm{in}}(X); \cal{G}) \to 0.$$ Under the generators $\lbrace \widetilde{\beta_{i}^{*}}\rbrace _{i=1}^{n}$ and $\lbrace \widetilde{\alpha_{i}^{*}}\rbrace _{i=1}^{n}$, the map $\widetilde{\partial}$ is represented by an $n\times n$ matrix $P$ which in turn gives a presentation of $A(G)$. 

An effective way to study $A(G)$ is to construct some values from the presentation matrix $P$ which do not depend on the choice of $P$.  
The $k$-\textit{th elementary ideal} of $P$ for $0\leq k\leq n$, which is denoted by $\epsilon_{k}(P)$, is the ideal of $\cal{G}$ generated by all $(n-k)\times (n-k)$ minors of $P$. It is known that $\epsilon_{k}(P)$ does not depend on the choice of the presentation matrix, and therefore is an invariant of $A(G)$. Since $\cal{G}$ is a unique factorization domain, we can define the $k$-\textit{th characteristic polynomial}, which is the greatest common divisor of all $(n-k)\times (n-k)$ minors of $P$. 

In this paper we are only interested in $\epsilon_{0}(P)$, which is principal and generated by the determinant of $P$. We call $\det(P)$ the \textit{Alexander invariant} of $G$. It is easy to see that $\det(P)$ modulo $\pm H_{1}(X; \mathbb{Z})$ does not depend on the choice of $P$ and therefore is an invariant of $G$. Unlike the case of $\theta_{n}$-graph, the determinant can be zero for some graphs.

\subsection{Fox calculus}
We show how to calculate the Alexander invariant of a graph $G$ using Fox calculus. 
When ${\partial}_{\mathrm{in}}(X)$ is disconnected, we follow the idea in \cite{MR699242} to construct a joint pair $(X', \partial_{\mathrm{in}}X')$ associated with $(X, \partial_{\mathrm{in}}X)$ as follows, where $\partial_{\mathrm{in}}X'$ is a connected subcomplex in $X'$.
Add a new $0$-cell $p$ to $X$ and add a new $1$-cell joining $p$ to a $0$-cell in each component of ${\partial}_{\mathrm{in}}(X)$. Then let $X'$ be the union of $X$, $p$ and these new $1$-cells, and let $\partial_{+}X'$ be the union of $\partial_{\mathrm{in}}X$, $p$ and these new $1$-cells. 

\begin{figure}
	\centering
		\includegraphics[width=0.3\textwidth]{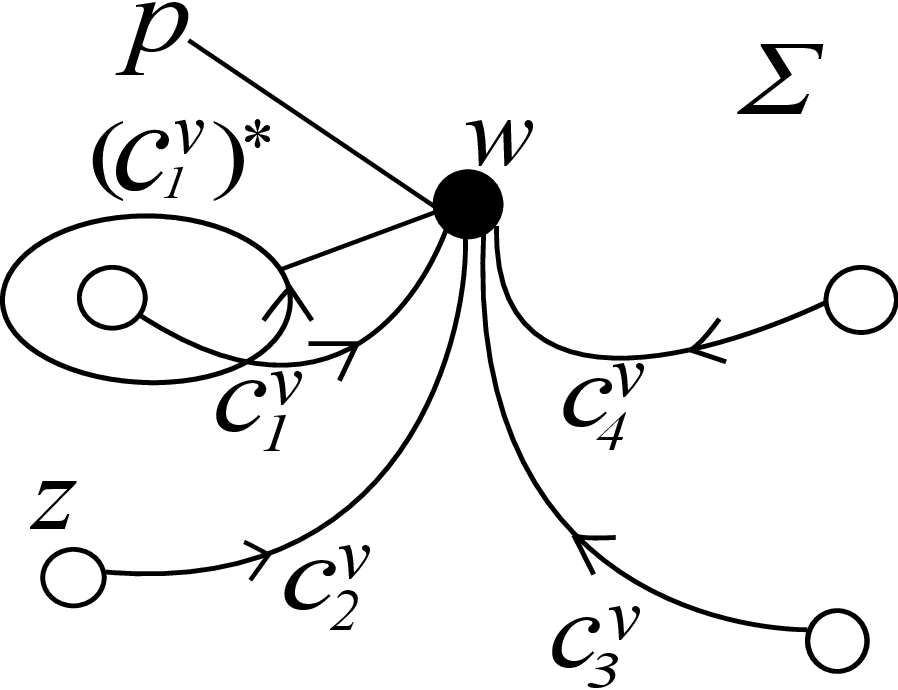}
		\caption{}
	\label{ff11}
\end{figure}

Now we discuss how to get a presentation for the fundamental group of $X'$ from a Heegaard diagram of $G$. Suppose $H_{G}=(\Sigma, \boldsymbol{\alpha}, \boldsymbol{\beta}, \boldsymbol{w}, \boldsymbol{z})$ is a Heegaard diagram of $G$ as before. We choose an orientation for each $\boldsymbol{\alpha}$- and $\boldsymbol{\beta}$-curve. As in Figure~\ref{ff11}, define the oriented arcs $c_{1}^{v}, c_{2}^{v}, \cdots, c_{d^{\mathrm{in}}(v)}^{v}$ for each vertex $v$ with indegree $d^{\mathrm{in}}(v)$, where $c_{i}^{v}$ is an oriented arc from a base point of type $\boldsymbol{z}$ to the base point $w_v$ around $v$. Note that $c_{i}^{v}$ for any $v\in V$ can be chosen to be disjoint from the $\boldsymbol{\alpha}$-curves. Let $(c_i^v)^*$ be the dual curve of $c_{i}^{v}$ in $\Sigma$ in the sense that it intersects with $c_{i}^{v}$ at a single point and is disjoint from any other curves of its own type. The curve $(c_i^v)^*$ is oriented so that its intersection point with $c_{i}^{v}$ has positive sign. We can construct a loop in $X'$ with base point $p$ by connecting $p$ to $(c_i^v)^*$. For simplicity, we still denote the loop by $(c_i^v)^*$.

Now we get a presentation of the fundamental group of $X'$ as follows.
$$\pi_{1}(X', p)=\Big \langle \lbrace \alpha_{i}^{*} \rbrace _{i=1}^{n}, \lbrace (c_i^v)^* \rbrace _{v\in V, 1\leq i \leq d^{\mathrm{in}}(v)}  \Big |  \lbrace \beta_{j}^{*} \rbrace _{j=1}^{n} \Big \rangle,$$ where $\alpha_{i}^{*}$ is a loop obtained by connecting $p$ to the 1-handle with belt circle $\alpha_{i}$. The relation $\beta_{j}^{*}$ can be easily read off from the Heegaard diagram $H_{G}$ as follows. Choose a generic point on $\beta_{j}$ as the start point and travel along $\beta_{j}$. When meeting an intersection point $x$ between $\alpha_{i}$ (resp. $c_{i}^{v}$) and $\beta_{j}$, we record $(\alpha_{i}^{*})^{{\mathrm {sign}}(x)}$ (resp. $[(c_{i}^{v})^{*}]^{{\mathrm {sign}}(x)}$).

We have the following proposition, the proof of which follows from the proof of Theorem 2.5 of \cite{MR699242}, although the context here is different from that of there.

\begin{prop}
\label{hempel}
Suppose that the component number of ${\partial}_{\mathrm{in}}(X)$ is $n$ (which is the number of vertices of $G$).
Then we have
\begin{enumerate}
\item  $\pi_{1}(X')\cong \pi_{1}(X)\ast F$, where $F$ is the free group of rank $n-1$, and $H_{1}(X, {\partial}_{\mathrm{in}}(X))\cong H_{1}(X', {\partial}_{\mathrm{in}}(X'))$;
\item For a Heegaard diagram $H_G$ of $G$, the matrix $ ([\dfrac{\partial \beta_{j}^{*}}{\partial \alpha_{i}^{*}}])_{i, j=1}^{n}$ represents the map $\widetilde{\partial}: C_{2}(\widetilde{X}, \widetilde{\partial}_{\mathrm{in}}(X); \cal{G}) \to C_{1}(\widetilde{X}, \widetilde{\partial}_{\mathrm{in}}(X); \cal{G}) $ discussed in Section 5.2, and therefore is a presentation matrix for the Alexander module $A(G)$, where $\dfrac{\partial \beta_{j}^{*}}{\partial \alpha_{i}^{*}}$ is calculated by applying the Fox calculus and $[\cdot]$ is the composition of the Hurewicz map from $\pi_{1}(X', p)$ to $H_{1}(X'; \mathbb{Z})$ and the projection from $H_{1}(X'; \mathbb{Z})$ to $H_{1}(X; \mathbb{Z})$.
\end{enumerate}
\end{prop}

On the other hand, following the proof of \cite[Prop. 4.2]{MR2805998}, it is easy to see the following theorem, which states that the Alexander invariant coincides with the Euler Characteristic of the Heegaard Floer complex of $G$. 

\begin{theo}[Prop. 4.2 \cite{MR2805998}] We have
$\det([\dfrac{\partial \beta_{j}^{*}}{\partial \alpha_{i}^{*}}])_{i, j=1}^{n}\doteq \chi(\widehat{\mathrm{CFG}}(\Sigma, \boldsymbol{\alpha}, \boldsymbol{\beta}, \boldsymbol{w}, \boldsymbol{z}))$.
\end{theo}

\subsection{State sum formula}
Given a graph diagram of $G$, we provide a diagrammatical formula for the Alexander invariant in terms of state sum, which is a generalization of Kauffman's state sum formula for the Alexander polynomial of a link.

\begin{figure}
\begin{tikzpicture}[baseline=-0.65ex, thick]
\draw (0, 1) to (1.75, 1) ;
\draw (2.25, 1) [->] to (4, 1) ;
\draw (2, -0.5) [->] to (2, 2.5) ;
\draw (2.35, 2.55) node {$e$};
\draw (2.5, 1.35) node {\large $1$};
\draw (2.5, 0.65) node {\large $1$};
\draw (1.6, 1.35) node {\large $t_e$};
\draw (1.6, 0.65) node {\large $t_e$};
\end{tikzpicture}
\quad\quad\quad
\begin{tikzpicture}[baseline=-0.65ex, thick]
\draw (1, 2) to [out=270,in=180] (2.5, 1.3);
\draw (4, 2) to [out=270,in=0] (2.5, 1.3);
\draw (4, 2) to [out=90,in=0] (2.5, 2.7);
\draw  (1, 2) to [out=90,in=180] (2.5, 2.7);
\draw (2.5, 1.3) [-<-] to (2.5, 0);
\draw (2.75, 0.2) node { $e$};
\draw (2.85, 1) node {\large $1$};
\draw (2.15, 1) node {\large $t_e$};
\draw (2.5, 1.65) node {\large $t_e-1$};
\end{tikzpicture}\\
\begin{tikzpicture}[baseline=-0.65ex, thick]
\draw (0, 1) to (1.75, 1) ;
\draw (2.25, 1) [->] to (4, 1) ;
\draw (2, -0.5) [->] to (2, 2.5) ;
\draw (2.35, 2.55) node {$e$};
\draw (2.5, 1.35) node {\large $-1$};
\draw (2.5, 0.65) node {\large $1$};
\draw (1.6, 1.35) node {\large $1$};
\draw (1.5, 0.65) node {\large $-1$};
\end{tikzpicture}
\quad\quad\quad
\begin{tikzpicture}[baseline=-0.65ex, thick]
\draw (1, 2) to [out=270,in=180] (2.5, 1.3);
\draw (4, 2) to [out=270,in=0] (2.5, 1.3);
\draw (4, 2) to [out=90,in=0] (2.5, 2.7);
\draw  (1, 2) to [out=90,in=180] (2.5, 2.7);
\draw (2.5, 1.3) [-<-] to (2.5, 0);
\draw (2.75, 0.2) node { $e$};
\draw (2.85, 1) node {\large $1$};
\draw (2.05, 1) node {\large $-1$};
\draw (2.5, 1.65) node {\large $1$};
\end{tikzpicture}
	\caption{Around each crossing $C_{j}$, the local contributions $A_{C_j}^{s(C_{j})}$ (top) and $m_{C_j}^{s(C_{j})}$ (bottom) are determined by the rules above, where $t_e$ is the homology class of the oriented meridian of the edge $e$.}
\label{e4}
\end{figure}
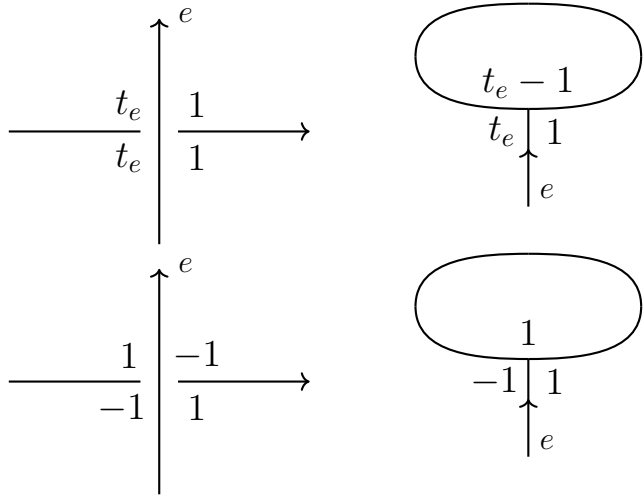

Suppose $G$ is a balanced bipartite graph in $S^{3}$, and $D$ is a connected graph diagram of $G$ in $S^2$. We define an element $\Delta_{D}\in \cal{G}$ from $D$ and prove that modulo $\pm H_{1}(X; \mathbb{Z})$ it coincides with the Alexander invariant of $G$. We finish the definition in the following steps.
\begin{enumerate}
\item The diagram $D$ separates $S^{2}$ into several regions. Choose a vertex $u\in V$. Asterisk the regions around $u$ which are adjacent to the edges pointing to $u$.
\item For any vertex $v\in V\backslash\lbrace u\rbrace$, place a small circle in $S^2$ centered at $v$. We call the disk bounded by the circle the \textit{circle region around $v$}. We call the intersection points of the circle with the edges pointing to $v$ the \textit{crossings around $v$}. There are totally $d^{\mathrm{in}}(v)$ crossings around $v$ where $d^{\mathrm{in}}(v)$ is the indegree of $v$.
\item Define $\mathrm{Cr}(D)$ to be the union of the double point crossings coming from $D$ and the crossings around each vertex created in Step (ii). Define $\mathrm{Re(D)}$ to be the union of the unasterisked regions and the circle regions.
\item If $\vert \mathrm{Cr}(D)\vert \neq \vert \mathrm{Re}(D)\vert$, which is the case in Lemma \ref{crre}, let $\Delta_{D}=0$. 
\item Now we suppose $p=\vert \mathrm{Cr}(D)\vert = \vert \mathrm{Re}(D)\vert$. Each double point is locally adjacent to four corners (which not necessarily belong to distinct regions), and each crossing around a vertex is locally adjacent to three corners. A \textit{state} of $D$ is a bijective map $s: \mathrm{Cr}(D) \to \mathrm{Re}(D)$ sending each crossing to one of its adjacent corners. 
\item Assign an order to the crossings (resp. regions) in $\mathrm{Cr}(D)$ (resp. $\mathrm{Re}(D)$). For a state $s$, define $\mathrm{sign}(s)$ to be the sign of $s$ as a permutation. 
Let $$A(s):=\prod_{j=1}^{p}A_{C_j}^{s(C_{j})}\in H_{1}(X; \mathbb{Z}), \text{ and } m(s):=\prod_{j=1}^{p}m_{C_j}^{s(C_{j})}\in \{1, -1\},$$ where $A_{C_j}^{s(C_{j})}$ and $m_{C_j}^{s(C_{j})}$ are defined by the rules in Fig.~\ref{e4}.
\item Define 
\begin{equation}
\label{tau}
\Delta_{D}:=\sum_{s: \text{ state }}\mathrm{sign}(s)\cdot m(s)\cdot A(s) \in \mathbb{Z}H_{1}(X; \mathbb{Z}).
\end{equation}

\end{enumerate}

The following lemma shows when the situation in Step (iv) can occur.

\begin{lemma}
\label{crre}
$\vert \mathrm{Cr}(D)\vert \neq \vert \mathrm{Re}(D)\vert $ if and only if the number of asterisked regions around $u$ is less than $d^{\mathrm{in}}(u)+1$.
\end{lemma}
\begin{proof}
Since $D$ is a connected diagram, we can calculate the Euler characteristic of $S^{2}$ using $D$. The number of $0$-simplices $p$ is the sum of the number of crossings of type \diaCrossP or \diaCrossN and the number of vertices of $D$. The number of $1$-simplices $q$ is the sum of the number of edges of $D$ and twice the number of crossings of type \diaCrossP or \diaCrossN. The number of $2$-simplices $r$ is the number of regions separated by $D$. We have $p-q+r=2$.

On the other hand, $\vert \mathrm{Cr}(D)\vert$ is the sum of the number of crossings of type \diaCrossP or \diaCrossN and the number of edges of $D$ minus $d^{\mathrm{in}}(u)$. Let $r^{*}$ be the number of asterisked regions. Then we see that $\vert \mathrm{Re}(D)\vert$ is the sum of the number of vertices of $D$ and $r$ minus $r^{*}+1$.

From the relation $p-q+r=2$ we see that $\vert \mathrm{Cr}(D)\vert=\vert \mathrm{Re}(D)\vert+r^{*}-(d^{\mathrm{in}}(u)+1)$, which implies the lemma. 
\end{proof}

\begin{figure}
	\centering
		\includegraphics[width=0.35\textwidth]{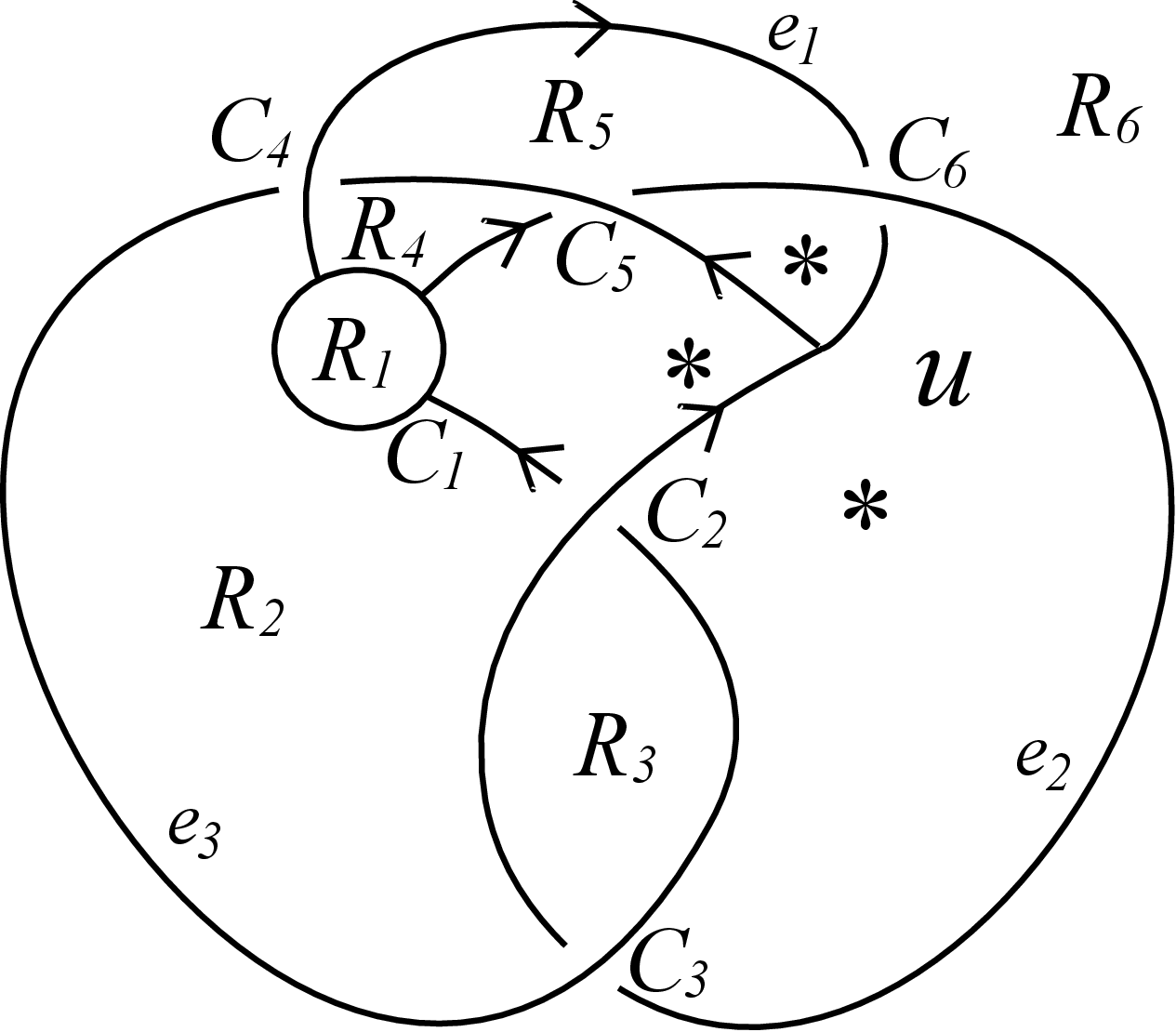}
	\caption{ Kinoshita's $\theta$-graph.}	
	\label{fig:ff3}
\end{figure}

\begin{ex}
For the Kinoshita's $\theta$-graph $D$ in Figure~\ref{fig:ff3}, we choose $u$ to be the vertex on the left. There are six crossings and six unasterisked regions. We call them $\{C_{1}, C_{2}, C_{3}, C_{4}, C_{5}, C_6\}$ and $\{R_{1}, R_{2}, R_{3}, R_{4}, R_{5}, R_6\}$ as in the figure. We calculate $\Delta_D$ as follows, where we represent a state by its underlying permutation and represent $t_{e_i}$ by $t_i$ for $i=1, 2, 3$. As a result, we have 
\begin{eqnarray*}
\Delta_D &\doteq& (t_3-1)t_2t_3(t_2-t_2t_3+t_1t_2t_3-1+t_3-t_1t_3+t_1)\\
&\doteq& (t_3-1)(t_2-t_2t_3+t_1t_2t_3-1+t_3-t_1t_3+t_1)\\
&\doteq& (t_1t_2-1)(t_1^2t_2^2-t_1t_2(t_1+t_2)+t_1t_2+t_1+t_2-1),
\end{eqnarray*}
where the third equality follows from the fact that $t_3=t_1t_2$.
\begin{table}[ht!]
  \begin{center}
\begin{tabular}{|c|c|c|c|}
  \hline
  Kauffman state $s$ & sign of $s$ & Value of $m(s)$ & Value of $A(s)$      \\
  \hline
  $123456$ & $1$ & $1$ &$(t_3-1)t_2^2 t_3$ \\
  $123546$ & $-1$ &$1$ &$(t_3-1)t_2^2 t_3^2$\\
  $123645$ & $1$ &$1$ &$(t_3-1)t_1t_2^2 t_3^2$\\
  $132456$ & $-1$ &$1$ &$(t_3-1)t_2 t_3$\\
  $132546$ & $1$ &$1$ &$(t_3-1)t_2 t_3^2$\\
  $132645$ & $-1$ &$1$ &$(t_3-1)t_1t_2 t_3^2$\\
  $136245$ & $1$ &$1$ &$(t_3-1)t_1t_2 t_3$\\
  \hline
\end{tabular}
\end{center}
\bigskip
  \label{table1}
  \end{table}
\end{ex}

\begin{rem}
The definition $\Delta_{D}$ depends on many choices: the vertex $u$, the order of crossings and that of regions, but these choices do not affect the value of $\Delta_{D}$ modulo $\pm H_{1}(X; \mathbb{Z})$, as we can see from the following theorem.
\end{rem}

In Example \ref{hdex} of Section 3.1 we showed how to construct a Heegaard diagram $H_{D}$ from a graph diagram $D$. We use it to prove the following theorem.
\begin{theo}
\label{eulertheo}
Modulo $\pm H_{1}(X; \mathbb{Z})$, the state sum $\Delta_{D}$ coincides with the Alexander invariant of $D$.
\end{theo}

\begin{proof}
Consider the Heegaard diagram $H_{D}$ in Example \ref{hdex}, where we choose to remove $\alpha_u$ and the $\boldsymbol{\alpha}$-curve on the right hand side of $u$, as shown in the following figure. 
\begin{figure}[!h]
	\centering
		\includegraphics[width=0.35\textwidth]{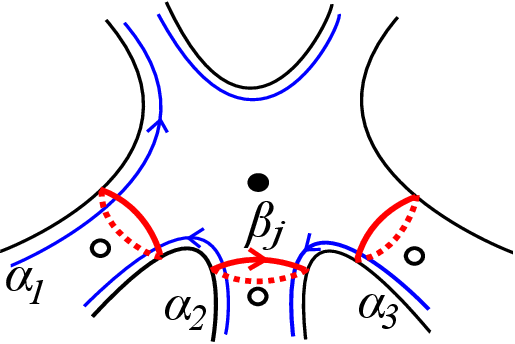}
\end{figure}

We prove that $\det ([\dfrac{\partial \beta_{j}^{*}}{\partial \alpha_{i}^{*}}])_{i,j=1}^{n}=\Delta_{D}$ modulo $\pm H_{1}(X; \mathbb{Z})$. For simplicity, we omit the symbol $*$ from $[\dfrac{\partial \beta_{j}^{*}}{\partial \alpha_{i}^{*}}]$ during the calculation.
We choose the counter-clockwise orientation for each $\boldsymbol{\alpha}$-curve. The orientations of $\boldsymbol{\beta}$-curves will be given in the following discussions. To calculate $[\dfrac{\partial \beta_{j}}{\partial \alpha_{i}}]$, depending on the position of $\beta_{j}$, we have three cases to consider.

\textbf{Case (i):} $\beta_{j}$ is around the vertex $u\in V$ (see the figure above). If the number of $\boldsymbol{\alpha}$-curves around $u$ that are adjacent to the edges pointing to $u$ is less than $d^{\mathrm{in}}(u)$,  it is easy to see that in this case the determinant of $([\dfrac{\partial \beta_{j}}{\partial \alpha_{i}}])_{i,j=1}^{n}$ is zero. This corresponds to the case of Step (iv) and Lemma \ref{crre} in the construction of $\Delta_{D}$. 

We now assume that the number of $\boldsymbol{\alpha}$-curves around $u$ that are adjacent to the edges pointing to $u$ is $d^{\mathrm{in}}(u)$. In this case $\beta_{j}=\alpha_{k+1}\alpha_{k}^{-1}(c_{k}^{u})^{-1}$ for some $1\leq k\leq d^{\mathrm{in}}(u)$ where we define $\alpha_{d^{\mathrm{in}}(u)+1}$ to be $1$. Note that $[\alpha_{k+1}]=[\alpha_{k}] t_{k}$, where $t_k$ is the homology class of the meridian for the edge adjacent to $\alpha_{k}$ and $\alpha_{k+1}$ around $u$. Therefore we have$$[\dfrac{\partial \beta_{j}}{\partial \alpha_{k}}]=-[\alpha_{k+1}\alpha_{k}^{-1}]=-t_{k}, \quad [\dfrac{\partial \beta_{j}}{\partial \alpha_{k+1}}]=1, $$ and $[\dfrac{\partial \beta_{j}}{\partial \alpha_{i}}]=0$ for $i\neq k, k+1$. For the $\beta_{j}$ which only intersects $\alpha_{d^{\mathrm{in}}(u)}$, we have $[\dfrac{\partial \beta_{j}}{\partial \alpha_{d^{\mathrm{in}}(u)}}]=-[\alpha_{d^{\mathrm{in}}(u)}^{-1}]$. If we let $I_{u}$ and $J_{u}$ denote the sets of indices of the $\boldsymbol{\alpha}$-curves around $u$ that are adjacent to the edges pointing to $u$ and  the $\boldsymbol{\beta}$-curves around $u$, respectively, we then have $\vert I_{u} \vert=\vert J_{u}\vert=d^{\mathrm{in}}(u)$. We see from the calculation above that modulo $\pm H_{1}(X; \mathbb{Z})$, the determinant of $([\dfrac{\partial \beta_{j}}{\partial \alpha_{i}}])_{i,j=1}^{n}$ is the same as that of $([\dfrac{\partial \beta_{j}}{\partial \alpha_{i}}])_{i, j=1, i\notin I_{u}, j\notin J_{u}}^{n}$. 

It is easy to see from the construction that $\{\alpha_{i}\}_{i=1, i\notin I_{u}}^{n}$ and $\{\beta_{j}\}_{j=1, j\notin J_{u}}^{n}$ correspond to the set $\mathrm{Re}(D)$ and the set $\mathrm{Cr}(D)$ in the definition of $\Delta_{D}$, respectively. We assume that $\alpha_{i}$ corresponds to $R_i\in \mathrm{Re}(D)$ and $\beta_{j}\in \mathrm{Cr}(D)$ corresponds to $C_j$. We show that the Leibniz formula for $\det([\dfrac{\partial \beta_{j}}{\partial \alpha_{i}}])_{i, j=1, i\notin I_{u}, j\notin J_{u}}^{n}$ coincides with the formula (\ref{tau}) of $\Delta_{D}$. 

\textbf{Case (ii): }$\beta_{j}$ corresponds to a double point (see the figure below). In this case, $\beta_{j}=\alpha_{a}\alpha_{b}^{-1}\alpha_{d}\alpha_{c}^{-1}$. The indices $a, b, c, d$ are not necessarily distinct. Let $\alpha_i=1$ if it does not exist for $i=a, b, c, d$. At any rate, the Fox calculus allows us to treat them as different indices in the calculation. 
\begin{figure}[ht!]
	\centering
		\includegraphics[width=0.25\textwidth]{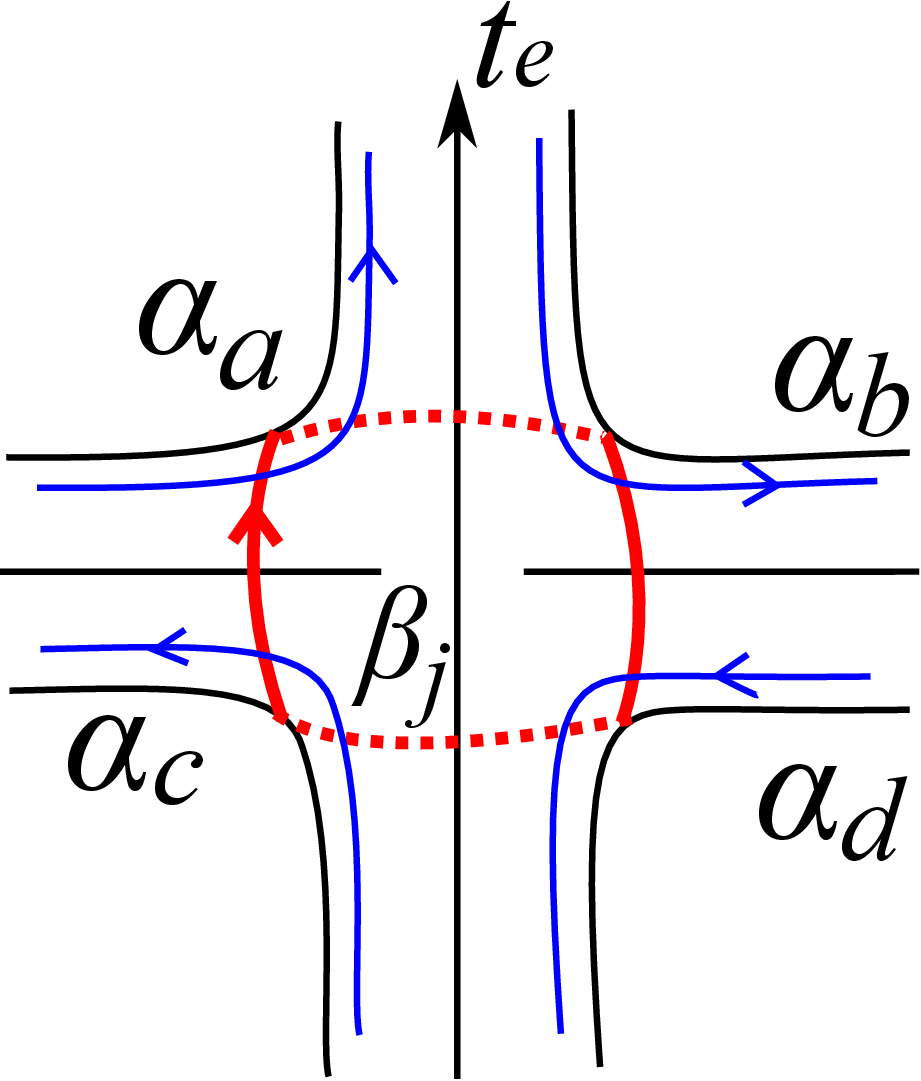}
\end{figure}

We see that
\begin{equation*}
\begin{aligned}
\dfrac{\partial \beta_{j}}{\partial \alpha_{a}} & =  1, \\\dfrac{\partial \beta_{j}}{\partial \alpha_{b}} &=-\alpha_{a}\alpha_{b}^{-1},
\end{aligned}\quad
\begin{aligned}
\dfrac{\partial \beta_{j}}{\partial \alpha_{c}}& = -\alpha_{a}\alpha_{b}^{-1}\alpha_{d}\alpha_{c}^{-1}, \\ \dfrac{\partial \beta_{j}}{\partial \alpha_{d}} &= \alpha_{a}\alpha_{b}^{-1}.
\end{aligned}
\end{equation*}
We have the relations $[\alpha_{b}]=[\alpha_{a}]t_{e}$ and $[\alpha_{d}]=[\alpha_{c}]t_{e}$. Therefore $$[\dfrac{\partial \beta_{j}}{\partial \alpha_{a}} ] =  1, [\dfrac{\partial \beta_{j}}{\partial \alpha_{b}}]=-t_{e}^{-1}, [\dfrac{\partial \beta_{j}}{\partial \alpha_{c}}]=-1, [\dfrac{\partial \beta_{j}}{\partial \alpha_{d}} ]=t_{e}^{-1},$$ and $[\dfrac{\partial \beta_{j}}{\partial \alpha_{i}}]=0$ for those $i\neq a, b, c, d$. Modulo $\pm H_{1}(X; \mathbb{Z})$, the determinant of the matrix $([\dfrac{\partial \beta_{j}}{\partial \alpha_{i}}])_{i, j=1, i\notin I_{u}, j\notin J_{u}}^{n}$ does not change if we let $$[\dfrac{\partial \beta_{j}}{\partial \alpha_{a}} ] =t_{e}, [\dfrac{\partial \beta_{j}}{\partial \alpha_{b}}]=-1, [\dfrac{\partial \beta_{j}}{\partial \alpha_{c}}]=-t_{e}, [\dfrac{\partial \beta_{j}}{\partial \alpha_{d}} ]=1,$$ and $[\dfrac{\partial \beta_{j}}{\partial \alpha_{i}}]=0$ for those $i\neq a, b, c, d$. We see that $[\dfrac{\partial \beta_{j}}{\partial \alpha_{i}}]=m_{C_j}^{R_i} \cdot A_{C_j}^{R_i}$ for $i=a, b, c, d$.

\textbf{Case (iii):} $\beta_{j}$ is around a vertex $v\in V-\lbrace u\rbrace$ (see the figure below). In this case $\beta_{j}$ intersects with at most three $\boldsymbol{\alpha}$-curves, one of which is $\alpha_v$ around $v$. Suppose that $\beta_{j}=\alpha_{k}^{-1}\alpha_{v}(c_{k}^{v})^{-1}\alpha_{v}^{-1}\alpha_{k+1}$ for some $1\leq k\leq d^{\mathrm{in}}(v)$. Let $\alpha_i=1$ if it does not exist for $i=k, k+1$.
\begin{figure}[ht!]
	\centering
		\includegraphics[width=0.4\textwidth]{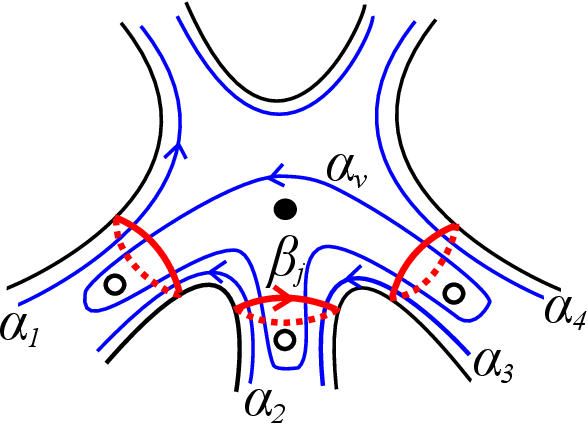}
\end{figure}

We have 
\begin{eqnarray*}
\dfrac{\partial \beta_{j}}{\partial \alpha_{k}} &=&-\alpha_{k}^{-1},\\
\dfrac{\partial \beta_{j}}{\partial \alpha_{k+1}}&=&\alpha_{k}^{-1}\alpha_{v}(c_{k}^{v})^{-1}\alpha_{v}^{-1},\\
\dfrac{\partial \beta_{j}}{\partial \alpha_{v}}&=&\alpha_{k}^{-1}(1-\alpha_{v}(c_{k}^{v})^{-1}\alpha_{v}^{-1}),
\end{eqnarray*}
and $[\dfrac{\partial \beta_{j}}{\partial \alpha_{i}}]=0$ for those $i\neq k, k+1, v$. Note that we have the relations $[\alpha_{k+1}]=[\alpha_{k}] t_{k}$ and $[c_{k}^{v}]=t_{k}$, where $t_k$ is the homology class of the meridian for the edge adjacent to $\alpha_{k}$ and $\alpha_{k+1}$ around $v$. Multiplying $[\alpha_{k}]t_{k}$ to the $j$-th column of the matrix $([\dfrac{\partial \beta_{j}}{\partial \alpha_{i}}])_{i, j=1, i\notin I_{u}, j\notin J_{u}}^{n}$, which will not change $\det ([\dfrac{\partial \beta_{j}}{\partial \alpha_{i}}])_{i, j=1, i\notin I_{u}, j\notin J_{u}}^{n}$ modulo $H_1 (X; \mathbb{Z})$, we get
$$
[\dfrac{\partial \beta_{j}}{\partial \alpha_{k}}] = -t_{k},
[\dfrac{\partial \beta_{j}}{\partial \alpha_{k+1}}] = 1, \text{ and }
[\dfrac{\partial \beta_{j}}{\partial \alpha_{v}}]= t_{k}-1,
$$
and $[\dfrac{\partial \beta_{j}}{\partial \alpha_{i}}]=0$ if $i\neq k, k+1, v$. It is easy to see that $[\dfrac{\partial \beta_{j}}{\partial \alpha_{i}}]=m_{C_j}^{R_i} \cdot A_{C_j}^{R_i}$ for $i=k, k+1, v$.

\end{proof}

\begin{rem}
(1) From the proof we see that the construction of $\Delta_{D}$ depends deeply on the fundamental group presentation, and more precisely on the Heegaard diagram $H_{D}$ and the orientations of the $\boldsymbol{\alpha}$- and the $\boldsymbol{\beta}$-curves.\\
(2) Roughly speaking, a state $s$ in the definition (\ref{tau}) of $\Delta_{D}$ corresponds to a generator $x$ of the Heeegaard Floer complex $\widehat{\mathrm{CFG}}(\Sigma, \boldsymbol{\alpha}, \boldsymbol{\beta},\boldsymbol{w}, \boldsymbol{z})$. In this sense, $\mathrm{sign}(\sigma)\prod_{i=1}^{n}\mathrm{sign}(x_{i})$ in (\ref{euler}) is $\mathrm{sign}(s)\cdot m(s)$ in (\ref{tau}), and $A(x)$ is $A(s)$.

\end{rem}

\bibliographystyle{siam}
\bibliography{bao}

\end{document}